\documentclass[11pt]{article}
\usepackage{amsmath,amsthm,latexsym,amsfonts, amssymb}
\usepackage{enumerate,epsfig}

\newcommand{\cB}{\mathcal{B}}

\newcommand{\nt}{{\lfloor nt\rfloor}}
\newcommand{\nns}{{\lfloor ns\rfloor}}

\newcommand{\N}{\mathbb N}
\newcommand{\R}{\mathbb R}
\newcommand{\Z}{\mathbb Z}
\newcommand{\Prob}{\mathbb P}
\newcommand{\E}{\mathbb E}
\newcommand{\eps}{\varepsilon}

\newcommand{\ra}{\rightarrow}
\newcommand{\lra}{\longrightarrow}
\newcommand{\ident}{{\mathchoice {\rm 1\mskip-4mu l} {\rm 1\mskip-4mu l}
{\rm 1\mskip-4.5mu l} {\rm 1\mskip-5mu l}}}

\newtheorem{theorem}{Theorem}[section]
\newtheorem{corollary}[theorem]{Corollary}
\newtheorem{lemma}[theorem]{Lemma}
\newtheorem{proposition}[theorem]{Proposition}
\newtheorem{definition}[theorem]{Definition}

\title{Asymptotic behaviour of the simple random walk\\
on the 2-dimensional comb
\thanks{Research partly supported by Italian 2004 PRIN project
``CAMPI ALEATORI''}}%[D. Bertacchi]
\author{Daniela Bertacchi\\
Universit{\`a} di Milano-Bicocca\\
 Dipartimento di Matematica e Applicazioni\\
 Via Cozzi 53, 20125 Milano, Italy\\
daniela.bertacchi\@@unimib.it }
\date{}

\begin{document}
\maketitle

\begin{abstract}
We analyze the differences between the horizontal and the vertical
component of the simple random walk on the 2-dimensional
comb. In particular
we evaluate by combinatorial methods
the asymptotic behaviour of the expected value of
the distance from the origin, 
the maximal deviation and the maximal span in $n$ steps,
proving that for all these quantities the order is $n^{1/4}$ for the 
horizontal projection and $n^{1/2}$ for the vertical one (the exact
constants are determined).
Then we rescale the two projections of the random walk dividing
by $n^{1/4}$ and $n^{1/2}$ the horizontal and vertical ones, 
respectively. The limit process is obtained.
As a corollary of the estimate of the expected value
of the maximal deviation, 
the walk dimension is determined, showing that the Einstein relation
between the fractal, spectral and walk dimensions does not hold on 
the comb. 
%+ 60F05 weak limits, clt; 60C05 combinatorial probability
\end{abstract}

%\maketitle

%%%%%%%%%%%%%%%%%%%%%%%%%%%%%%%%%%%%%%%%%
%%%	INTRODUCTION		      %%%
%%%%%%%%%%%%%%%%%%%%%%%%%%%%%%%%%%%%%%%%%

\section{Introduction and main results}
\label{sec:intro}
\setcounter{equation}{0}

\begin{figure}
\centerline{
\epsfig{figure=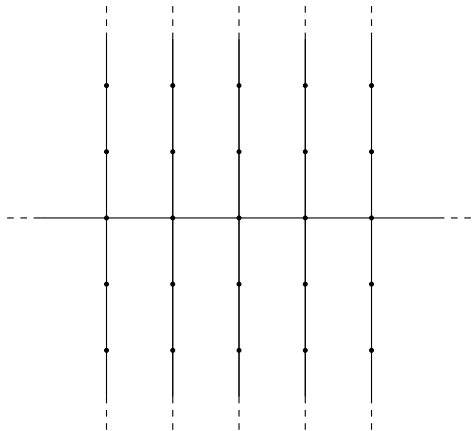,height=6cm}}
\caption{The 2-dimensional comb.}
\label{fig:2comb}
\end{figure}

The 2-dimensional comb $\mathbf C_2$
is maybe the simplest example of inhomogeneous graph.
It is obtained from $\Z^2$ by removing all horizontal edges
off the $x$-axis (see Figure~\ref{fig:2comb}).
Many features of the simple random walk on this graph has
been matter of former investigations. Local limit theorems
were first obtained by Weiss and Havlin \cite{Weiss-Havlin}
and then extended to higher dimensions by Gerl \cite{Gerl2}
and Cassi and Regina \cite{Cassi-Regina}. More recently,
Krishnapur and Peres \cite{cf:KP} have shown that on $\mathbf C_2$
 two independent walkers meet only finitely many times
almost surely. This result, together with the space-time
asymptotic estimates obtained in \cite{cf:BZ} for the $n$-step
transition probabilities, suggests that the walker spends most
of the time on some tooth of the comb, that is moving along
the vertical direction.
Indeed in \cite[Section 10]{cf:BZ}, it has been remarked that, 
if $k/n$ goes to zero with a certain speed, then
$p^{(2n)}\big((2k,0), (0,0)\big)/p^{(2n)}\big((0,2k), (0,0)\big)
\stackrel{n\to\infty}\lra0$.

Moreover the results in \cite{cf:BZ} imply that there are
no sub-Gaussian estimate of the transition probabilities on
$\mathbf C_2$. Such estimates have been found on many graphs:
by Jones \cite{Jo} on the $2$-dimensional Sierpi\'nski graph,
by Barlow and Bass \cite{Ba-Ba2} on the graphical Sierpi\'nski carpet
and on rather general graphs by Grigor'yan and Telcs
(\cite{Grig,Telcs}). These estimates involve three exponents which are
usually associated to infinite graphs: 
the {\it spectral dimension} $\delta_s$ (which is by definition twice
the exponent of $n^{-1}$ in local limit theorems), the
{\it fractal dimension} $\delta_f$ (which is the growth exponent)
 and  the {\it walk dimension} $\delta_w$. 
These dimensions are in typical cases linked
by the so-called {\it Einstein relation}:
$\delta_s\delta_w=2\delta_f$
(see Telcs \cite{Te1,Te2}).
The first two dimensions are known for $\mathbf C_2$: $\delta_s=3/2$
and $\delta_f=2$. In this paper we compute $\delta_w=2$, thus showing 
that the relation does not hold for this graph.

In order to point out the different behaviour of the random walk
along the two directions we analyze the asymptotic behaviour of
the expected value of the distance from the origin
reached by the walker in 
$n$ steps. Different concepts of distance are considered:
position after $n$ steps, maximal deviation and maximal span.
Section~\ref{sec:prelim} is devoted to these definitions, and to the
necessary preliminaries, such as the definitions of random
walk on a graph and of generating function. The
expression of the generating function of the transition probabilities
of the simple random walk on $\mathbf C_2$ is recalled.

In Sections~\ref{sec:distances}, \ref{sec:maxdev} and \ref{sec:maxspan}
we prove the asymptotic estimates of the expected value of the
distance of the walk from the origin after $n$ steps, of its
maximal deviation from the origin and of its maximal span respectively.
The proofs are based on a Darboux type transfer theorem:
we refer to \cite[Corollary 2]{cf:FO}, but one may also refer to
\cite{cf:Bender} and to the Hardy-Littlewood-Karamata theorem
(see for instance \cite{cf:Feller2}). 
This theorem (as far as we are concerned) claims that if
$$
F(z):=\sum_{n=0}^\infty a_nz^n \stackrel{z\ra 1^-}\sim \frac C{(1-z)^\alpha},
\qquad \alpha\not\in\{0,-1,-2,\ldots\},
$$
and $F(z)$ is analytic in some domain, with the exception of $z=1$,  then
$$
a_n\stackrel{n\to\infty}\sim \frac C{\Gamma(\alpha)}\, n^{\alpha-1}.
$$
The aim of our computation is then to determine an explicit expression
of the generating functions of the sequence of expected values of the
random variables we are interested in.
This is done employing the combinatorial methods used by Panny and Prodinger
in \cite{cf:PP}. We also refer to that paper for a comparison between our
results and the analogous results for the simple random walk on $\Z$.
Theorems~\ref{th:distx}, \ref{th:disty}, \ref{th:mdx}, \ref{th:mdy}, 
\ref{th:msx} and \ref{th:msy}
prove that the expected values of the distances along the horizontal
direction are all of order $n^{1/4}$, and along the vertical direction
they are of order $n^{1/2}$ (the exact constants are determined).

Since the notions of maximal deviation in $n$ steps and 
of first exit time from a $k$-ball are closely related,
we determine the walk dimension of $\mathbf C_2$ in 
Section~\ref{sec:maxdev}.

In Section~\ref{sec:scaling} we deal with the limit of the process
obtained dividing by $n^{1/4}$ and by $n^{1/2}$ respectively the continuous
time interpolation of the horizontal
and vertical projections of the position after $n$ steps. As one would
expect the limit of the vertical component is the Brownian motion, 
while the limit of the horizontal component is less obvious
(it is a Brownian motion indexed by the local time at 0 of the vertical
component).
This scaling limit is determined in Theorem~\ref{th:scaling}.

Finally, Section~\ref{sec:fr} is devoted to
 a discussion of the results, remarks and open questions.

%%%%%%%%%%%%%%%%%%%%%%%%%%%%%%%%%%%%%%%%%
%%%	PRELIMINARIES                 %%%
%%%%%%%%%%%%%%%%%%%%%%%%%%%%%%%%%%%%%%%%%

\section{Preliminaries}
\label{sec:prelim}

The simple random walk on a graph is a sequence of random variables
$\{S_n\}_{n\ge0}$, where $S_n$ represent the position of the walker 
at time $n$, such that if $x$ and $y$ are vertices which are
neighbours  then
$$
p(x,y):=\Prob\left(S_{n+1}=y|S_n=x\right)=\frac1{\mathrm{deg}(x)},
$$
where $\mathrm{deg}(x)$ is the number of neighbours of $x$,
otherwise $p(x,y)=0$.
In particular on $\mathbf C_2$
the non-zero transition probabilities $p(x,y)$ are
equal to 1/4 if $x$ is on the horizontal axis, 
and they are equal to 1/2 otherwise. %if $x$ is off the horizontal axis.

Given $x,y\in\mathbf C_2 $, let
$$
p^{(n)}(x,y):=\Prob\left(S_n=y|S_0=x\right),\qquad n\ge 0,
$$
be the $n$-step transition probability from $x$ to $y$.
Recall that the generating function of a sequence $\{a_n\}_{n\ge0}$
is the power series $\sum_{n\ge0}a_nz^n$; by definition the Green
function associated to the random walk on a graph $X$
is the family of generating
functions of the sequences $\{p^{(n)}(x,y)\}_{n\ge0}$, $x,y\in X$,
that is
$$
G(x,y|z)=\sum_{n\ge0}p^{(n)}(x,y)z^n.
$$
The Green function, with $x=(0,0)$, 
on $\mathbf C_2$ can be written explicitly as
(see \cite{cf:BZ})
$$
G\big((0,0),(k,l)|z\big)=
\begin{cases}
\frac12\, G(z)(F_1(z))^{|k|}(F_2(z))^{|l|},& \mathrm{if}\ l\neq0,\\
G(z)(F_1(z))^{|k|},&\mathrm{if }\ l=0,
\end{cases}
$$
where
\begin{equation*}
\begin{split}
 G(z) & =\frac{\sqrt2}{\sqrt{1-z^2+\sqrt{1-z^2}}};\\
F_1(z) & =\frac{1+\sqrt{1-z^2}-\sqrt 2
    \sqrt{1-z^2+\sqrt{1-z^2}}}{z};\\
F_2(z) & = \frac{1-\sqrt{1-z^2}}{z}.
\end{split}
\end{equation*}

We refer to \cite[Section 1.1]{cf:Woess} for more details on
the random walks on graphs, transition probabilities and generating
functions.

When we consider the walk up to time $n$, different concepts
of ``distance'' arise. 
We consider two equivalent norms on $\mathbf C_2$: for a given vertex
$(x,y)$ define
\begin{equation*}
\Vert (x,y)\Vert_1=|x|+|y|,\qquad
\Vert (x,y)\Vert_\infty=\max\{|x|,|y|\}.
\end{equation*}
Note that $\Vert \cdot\Vert_1$ is the usual distance on the graph.
In the following section we deal not only with the
asymptotic behaviour of $\E[|S_n^x|]$ and $\E[|S_n^y|]$ 
(we use the notation $S_n=(S_n^x,S_n^y)$), but also
with the asymptotics of the expected value of other random
variables, which represent the (horizontal and vertical)
maximal deviation and the span of the walk.
\begin{definition}
\begin{enumerate}[$a.$]
\item The maximal deviations in $n$ steps are defined as
\begin{equation*}
\begin{split}
D_n^x & :=\max\{|S_i^x|:0\le i\le n\};\\
D_n^y & :=\max\{|S_i^y|:0\le i\le n\}.\\
\end{split}
\end{equation*}
\item The maximal spans in $n$ steps are defined as
\begin{equation*}
\begin{split}
M_n^x & :=\max\{S_i^x-S_j^x:0\le i,j\le n\},\\
M_n^y & :=\max\{S_i^y-S_j^y:0\le i,j\le n\}.
\end{split}
\end{equation*}
\end{enumerate}
\end{definition}

%%%%%%%%%%%%%%%%%%%%%%%%%%%%%%%%%%%%%%%%%
%%%	DISTANCE    		      %%%
%%%%%%%%%%%%%%%%%%%%%%%%%%%%%%%%%%%%%%%%%

\section{Mean distance}
\label{sec:distances}

\begin{theorem}\label{th:distx}
\begin{equation*}
\E[|S_n^x|]\stackrel{n\to\infty}\sim\frac 1{2^{3/4}\Gamma(5/4)}\,n^{1/4}.
\end{equation*}
\end{theorem}

\begin{proof}
Since for $k\neq0$,
\begin{equation*}
\Prob \left(|S_n^x|=k\right) = 
2\sum_{l\in\Z} p^{(n)} \big((0,0),(k,l)\big),
\end{equation*}
it is clear that (exchanging the order of summation)
\begin{equation*}
\sum_{n=0}^\infty\E[|S_n^x|] z^n=
2\sum_{k=1}^\infty k\sum_{l\in\Z}
G\big((0,0),(k,l)|z\big).
\end{equation*}
 By elementary computation
one obtains
\begin{equation*}\begin{split}
\sum_{n=0}^\infty\E[|S_n^x|] z^n  &= 
G(z)\,\frac{F_1(z)}{(1-F_1(z))^2}\,\frac 1{1-F_2(z)}\\
&\stackrel{z\ra 1^-}\sim
\frac 1{2^{3/4}(1-z)^{5/4}}.
\end{split}
\end{equation*}
Thus, applying \cite[Corollary 2]{cf:FO} we obtain the thesis.
\end{proof}

\begin{theorem}\label{th:disty}

\begin{equation*}
\E[|S_n^y|]\stackrel{n\to\infty}\sim 
\sqrt{\frac 2\pi}\,
n^{1/2}.
\end{equation*}
\end{theorem}

\begin{proof}
As in the proof of the previous theorem, 
for $l\neq0$, we write 
\begin{equation*}
\Prob \left(|S_n^y|=l\right) = 
2\sum_{k\in\Z} p^{(n)} \big((0,0),(k,l)\big),
\end{equation*}
and 
\begin{equation*}
\sum_{n=0}^\infty\E[|S_n^y|] z^n =
2\sum_{l=1}^\infty l\sum_{k\in\Z}
G\big((0,0),(k,l)|z\big),
\end{equation*}
Thus
\begin{equation*}\begin{split}
\sum_{n=0}^\infty\E[|S_n^y|] z^n &=G(z)\,\frac{1+F_1(z)}{1-F_1(z)}
\, \frac{F_2(z)}{(1-F_2(z))^2}\\
&\stackrel{z\ra 1^-}\sim\frac 1{2^{1/2}(1-z)^{3/2}}.
\end{split}
\end{equation*}
Apply \cite[Corollary 2]{cf:FO}, recalling that $\Gamma(3/2)=\sqrt\pi/2$,
to conclude.
\end{proof}

\begin{corollary}
$\E[\Vert S_n\Vert_1]$ and 
$\E[\Vert S_n\Vert_\infty]$ are both asymptotic, as $n$ goes to
infinity, to $\sqrt{\frac 2\pi}\,
n^{1/2}$.
\end{corollary}

%%%%%%%%%%%%%%%%%%%%%%%%%%%%%%%%%%%%%%%%%
%%%	DEVIATION   		      %%%
%%%%%%%%%%%%%%%%%%%%%%%%%%%%%%%%%%%%%%%%%
\section{Mean maximal deviation and walk dimension}
\label{sec:maxdev}

\subsection{Maximal horizontal deviation}
\label{sec:maxdevx}

In order to compute the generating function of $\{\E[D_n^x]\}_{n\ge0}$,
we first need an expression for another generating function.

\begin{lemma}\label{th:pdnx}
Let $h\in\N\cup\{0\}$, $l\in\{0,\ldots,h\}$.
The generating function of the sequence
$\left\{\Prob \left(D_n^x\le h, S_n^x=l\right)\right\}_{n\ge0}$ is 
\begin{equation*}
\psi_{h,l}\left(\frac{1-\sqrt{1-z^2}}{2z}
\right)\cdot
\frac{2(1-\sqrt{1-z^2})}{z(\sqrt{1-z^2}-1+z)},
\end{equation*}
where $\psi_{h,l}(z)$ is the generating function of
the number of paths on $\Z$ of length $n$, from $0$ to $l$ with
maximal deviation less or equal to $h$.
\end{lemma}

\begin{proof}
Note that the paths we are interested in have no bound on the vertical
excursions. Thus we may decompose the walk into its horizontal 
and vertical components, and consider each horizontal step as a 
vertical excursion (whose length might be zero) coming back to the
origin plus a step along the horizontal direction.

Keeping this decomposition in mind, it is clear that the generating function 
of the sequence 
$\{\Prob\left(D_n^x\le h, S_n=(l,0)\right)\}_{n\ge0}$  is
\begin{equation*}
\psi_{h,l}\left(\frac{z\widetilde G(0,0|z)}4\right),
\end{equation*}
where $\widetilde G(0,j|z)$ is the generating
function of the probabilities of the $n$-step excursions
along one single tooth of the comb
(that is paths in Figure~\ref{fig:vertcomp} 
which do not use the loop at zero), from $(x,0)$ to $(x,j)$.
\begin{figure}
\centerline{\epsfig{file=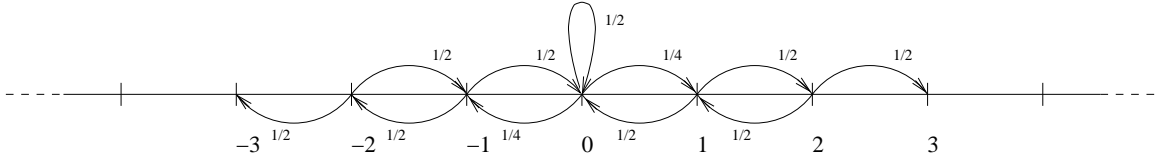,height=2.3cm}}
\caption{
The vertical component of the simple random walk on $\mathbf C_2$.
}
\label{fig:vertcomp}
\end{figure} 

Moreover we must admit  a final excursion from $(l,0)$ to 
$(l,j)$ for some $j\in\Z$, that is we must multiply by
\begin{equation*}
E(z):=
\widetilde G(0,0|z)+2\sum_{j=1}^\infty \widetilde G(0,j|z)
,
\end{equation*}
and the generating function of
$\left\{\Prob \left(D_n^x\le h, S_n^x=l\right)\right\}_{n\ge0}$ is 
\begin{equation*}
\psi_{h,l}\left(\frac{z\widetilde G(0,0|z)}{4}
\right)\cdot
E(z).
\end{equation*}
Using \cite[Lemma 1.13]{cf:Woess} and reversibility 
(see \cite[Section 1.2.A]{cf:Woess}),
it is not difficult to compute
\begin{equation*}
\begin{split}
\widetilde G(0,0|z) & = \frac{2\left(1-\sqrt{1-z^2}\right)}{z^2},\\
\widetilde G(0,j|z) & = \frac12\,\widetilde G(0,0|z)\left(
\frac{1-\sqrt{1-z^2}}z\right)^j,\quad j\neq0.
\end{split}
\end{equation*}
The thesis is obtained noting that
\begin{equation*}
\begin{split}
\frac{z\widetilde G(0,0|z)}{4} & =
\frac{1-\sqrt{1-z^2}}{2z},\\
E(z) &=\frac{2(1-\sqrt{1-z^2})}{z(\sqrt{1-z^2}-1+z)}.
\end{split}
\end{equation*}
\end{proof}

\begin{proposition}\label{th:genfundxv}
\begin{equation*}
\sum_{n=0}^\infty \E[D^x_n]z^n=
2\,\frac{1+6v^2+v^4}{(1-v)^4}\sum_{h\ge1}\frac{v^h}{1+v^{2h}},
\end{equation*}
where $v$ is such that 
\begin{equation}\label{eq:v&z}
\frac{v}{1+v^2}=\frac{1-\sqrt{1-z^2}}{2z}.
\end{equation}
\end{proposition}

\begin{proof}
Let $\psi_{h,l}(z)$ be as in Lemma~\ref{th:pdnx}, 
and put $\psi_h(z)=\sum_{|l|\le h}\psi_{h,l}(z)$.
Then the generating function of 
$\{\Prob \left(D_n^x\le h\right)\}_{n\ge0}$ is 
\begin{equation*}
H_h(z):=\psi_{h}\left(\frac{1-\sqrt{1-z^2}}{2z}\right)\cdot
\frac{2\left(1-\sqrt{1-z^2}\right)}{z(\sqrt{1-z^2}-1+z)}.
\end{equation*}
Thus we may write the generating function of $\E[D_n^x]$ as
\begin{equation}\label{eq:EDx}
\sum_{n=0}^\infty 
\left(\sum_{h=0}^\infty\Prob\left(D_n^x>h\right)\right)z^n
= \sum_{h=0}^\infty 
\left(\frac 1{1-z}-H_h(z)
\right).
\end{equation}
An explicit expression for $\psi_h$ has
been determined by Panny and Prodinger in \cite[Theorem 2.2]{cf:PP}:
\begin{equation*}
\psi_{h}(w)=\frac{(1+v^2)(1-v^{h+1})^2}{(1-v)^2(1+v^{2h+2})},
\end{equation*}
where $w=v/(1+v^2)$.
By the definition of $H_h(z)$ it is clear that the relation between
$z$ and $v$ is set by equation~\eqref{eq:v&z}.
Then it is just a matter of computation to obtain
\begin{equation*}\begin{split}
H_h(z)&=\frac{(1+6v^2+v^4)(1-v^{h+1})^2}{(1-v)^4(1+v^{2(h+1)})},\\
\frac1{1-z}& =\frac{1+6v^2+v^4}{(1-v)^4},
\end{split}\end{equation*}
whence, substituting in \eqref{eq:EDx}, the thesis.
\end{proof}

\begin{lemma}\label{th:vasint}
\begin{equation*}
\sum_{h\ge1}\frac{v^h}{1+v^{2h}}
\stackrel{v\ra 1^-}\sim\frac{\pi}{4(1-v)}.
\end{equation*}
\end{lemma}

\begin{proof}
The proof is quite standard, we report it here for completeness.
Put $v=e^{-t}$, $g(w)=e^{-w}/(1+e^{-2w})$ and consider
\begin{equation*}
f(t):=\sum_{h\ge1}\frac{e^{-ht}}{1+e^{-2ht}}=\sum_{h\ge1}g(ht).
\end{equation*}
The Mellin transform of $f$ is:
\begin{equation*}\begin{split}
f^*(s)& =\int_0^\infty\sum_{h\ge1}g(ht)\,t^{s-1} dt\\
& =
\zeta(s)\,\int_0^\infty \sum_{\lambda\ge0} (-1)^\lambda
e^{-(2\lambda+1)w}w^{s-1} dw,
\end{split}\end{equation*}
where $\zeta(s)=\sum_{h\ge1}h^{-s}$ is the Riemann zeta function.
The knowledge of the behaviour of $f^*(s)$ in a neighbourhood of 1, will
give us the behaviour of $f(t)$ in a neighbourhood of 0.
Substitute $y=(2\lambda+1)w$ to obtain
\begin{equation*}
f^*(s) = \zeta(s)\kappa(s)\Gamma(s),
\end{equation*}
where $\kappa(s)=\sum_{\lambda\ge0}\frac{(-1)^\lambda}{(1+2\lambda)^s}$,
and $\Gamma(s)= \int_0^\infty e^{-y}y^{s-1}dy$ is the gamma function.
Since as $s\ra 1^+$
\begin{equation*}
\begin{split}
\zeta(s) & =\frac 1{s-1}+O(1),\\
\Gamma(s)& =1+ O(s-1),\\
\end{split}
\end{equation*}
we are left with the computation of the asymptotic behaviour of
$\kappa(s)$. We may write
\begin{equation*}
\kappa(s) 
= \frac 1{4^s}\left(\zeta(s,1/4)- \zeta(s,3/4)\right),
\end{equation*}
where $\zeta(s,a)=\sum_{h\ge0}(a+h)^{-s}$ is the Hurwitz zeta function.
Thus using the expansion of $\zeta(s,a)$ for $s$ close to $1$
(see \cite[Formula 13.21]{cf:WW}) we obtain 
\begin{equation*}
\kappa(s) 
=\frac\pi4+O(s-1).
\end{equation*}
Hence we get
\begin{equation*}
f^*(s)\stackrel{s\ra 1^+}{\sim} \frac{\pi}{4(s-1)}.
\end{equation*}
Applying \cite[Theorem 1, p.115]{cf:Doetsch},
\begin{equation*}
f(t)\stackrel{t\ra0^+}{\sim} \frac\pi {4t},
\end{equation*}
which, substituting $t=-\log(1-(1-v))$, gives the thesis.
\end{proof}

\begin{theorem}\label{th:mdx}
\begin{equation*}
\E[D^x_n]\stackrel{n\to\infty}\sim
\frac{2^{-7/4}\pi}{\Gamma(5/4)}\,n^{1/4}.
\end{equation*}
\end{theorem}

\begin{proof}
By Proposition~\ref{th:genfundxv} and Lemma~\ref{th:vasint},
it is clear that
\begin{equation}\label{eq:dnxv}
\sum_{n\ge0}\E[D_n^x]z^n\stackrel{v\ra 1^-}{\sim}
\frac{4\pi}{(1-v)^5}.
\end{equation}
Choosing the solution  $v$ of equation \eqref{eq:v&z}
which is smaller than 1 when $z$ is smaller than 1,
and substituting it in \eqref{eq:dnxv} we obtain 
\begin{equation*}
\sum_{n=0}^\infty \E[D^x_n]z^n
\stackrel{z\ra 1^-}{\sim}
\frac {2^{-7/4}\pi}{(1-z)^{5/4}},
\end{equation*}
and (by \cite[Corollary 2]{cf:FO}) the thesis.
\end{proof}

\subsection{Maximal vertical deviation}
\label{sec:maxdevy}

We recall a result which can be found in 
\cite[Theorem 2.1]{cf:PP}, which is useful in the sequel.

\begin{lemma}\label{th:ai}
Let $\mathbf A_i(z/2)$ be the $(i+1)\times (i+1)$ matrix
\begin{equation*}
\left[\begin{array}{ccccc}
1     &  -z/2  & 0    & \cdots & 0\\
-z/2  &  1     & -z/2     & \      & \vdots\\
0     & \ddots & \ddots   & \ddots & 0 \\
\vdots& \      & \ddots   & \ddots & -z/2 \\
0     & \cdots & 0   & -z/2   & 1
\end{array}\right]
\end{equation*}
and let $a_i(z/2)$ be its determinant. Then
\begin{equation}\label{eq:ai}
a_i(z/2)=\frac{1-v^{2i+4}}{(1-v^2)(1+v^2)^{i+1}},
\end{equation}
where
\begin{equation}\label{eq:v&z2}
\frac{v}{1+v^2}=\frac z2.
\end{equation}
\end{lemma}

For the sake of simplicity, in the sequel we write $\mathbf A_i$ and
$a_i$ instead of $\mathbf A_i(z/2)$ and $a_i(z/2)$ respectively.

\begin{lemma}
Let $h\in\N\cup\{0\}$, $l\in\{0,\ldots,h\}$.
The generating function of the sequence
$\{\Prob\left(D^y_n\le h,S_n^y=l\right)\}_{n\ge0}$ is 
\begin{equation}\label{eq:psihat1}
\widehat\psi_{h,l}(z)=\frac{(z/2)^la_{h-l-1}}{(1-z/2)a_{h-1}-z^2
a_{h-2}/4}.
\end{equation}
\end{lemma}

\begin{proof}
Consider the absolute value of the vertical projection of the random walk 
as the random walk on the non negative integers
with one-step transition probabilites described in
Figure~\ref{fig:absvert}.
\begin{figure}
\centerline{\epsfig{file=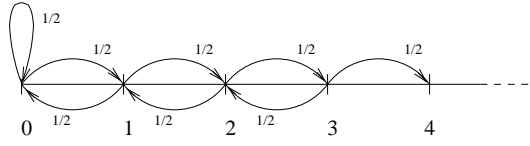,height=2.3cm}}
\caption{
The absolute value of the vertical component of the simple random walk on $\mathbf C_2$.}
\label{fig:absvert}
\end{figure} 
Note that $\widehat\psi_{h,l}(z)$, $l=0,\ldots,h$,
is determined by the linear system:
\begin{equation*}
\left[
\begin{array}{c} 
  \begin{array}{cccccc} %prima riga + seconda vuota
  (1-z/2) & -z/2 & 0       & \cdots  & 0 & \ \\
  \end{array}\\
\begin{array}{cc} % sotto due colonne, di cui una matrice
\begin{array}{c}
-z/2\\
0\\
\vdots\\
0\\
\end{array}
& \left[\begin{array}{lcccr}
\       &\       & \                & \ & \ \\
\       &\       & \                & \ & \ \\
\       &\       & \mathbf A_{h-1}  & \ & \ \\
\       &\       &\                 & \ & \ \\
\       &\       & \                & \ & \ \\
\end{array}\right]
\end{array}
\end{array}
\right]
\left[\begin{array}{c}
\widehat\psi_{h,0}(z) \\ 
\vdots \\\vdots \\
\widehat\psi_{h,h}(z) \\
\end{array}\right]
=
\left[\begin{array}{c}
1 \\ 0 \\ \vdots \\
0 \\
\end{array}\right].
\end{equation*}
Using Cramer's rule (compare with \cite[Theorem 2.1]{cf:PP}), 
we obtain the thesis.
\end{proof}

\begin{lemma}
The generating function of $\{\Prob[D^y_n\le h]\}_{n\ge0}$ is 
(written as a function of $v$)
\begin{equation*}
\widehat\psi_h(z)=\frac{(1+v^2)(1-v^{h+1})(1-v^{h+2})}{(1-v)^2(1+v^{2h+3})}.
\end{equation*}
\end{lemma}

\begin{proof}
From \eqref{eq:psihat1} and \eqref{eq:ai} we get
\begin{equation*}
\widehat\psi_{h,l}(z)=\frac{v^l(1+v^2)(1-v^{2h-2l+2})}{(1-v)(1+v^{2h+3})},
\end{equation*}
whence, summing $l$ from $0$ to $h$, the proof is complete.
\end{proof}

Proposition~\ref{th:genfundyv}, Lemma~\ref{th:vasint2} and
Theorem~\ref{th:mdy} are
the analogs of Proposition~\ref{th:genfundxv}, Lemma \ref{th:vasint}
and Theorem~\ref{th:mdx}
respectively. Therefore, their proofs are omitted.

\begin{proposition}\label{th:genfundyv}
Let $v$ be such that \eqref{eq:v&z2} holds.
Then
\begin{equation}
\sum_{n=0}^\infty \E[D_n^y]z^n=
\frac{(1+v^2)(1+v)}{(1-v)^2}\sum_{h\ge1}
\frac{v^{h}}{(1+v^{2h+1})}.
\end{equation}
\end{proposition}

\begin{lemma}\label{th:vasint2}
\begin{equation*}
\sum_{h\ge1}\frac{v^h}{1+v^{2h+1}}
\stackrel{v\ra 1^-}\sim\frac{\pi}{4(1-v)}.
\end{equation*}
\end{lemma}

\begin{theorem}\label{th:mdy}
\begin{equation*}
\E[D^y_n]\stackrel{n\to\infty}\sim
\sqrt{\frac\pi2}\, n^{1/2}.
\end{equation*}
\end{theorem}

\begin{corollary}
Both $\E[\max_{0\le i\le n}\Vert S_i\Vert_1]$
and $\E[\max_{0\le i\le n}\Vert S_i\Vert_\infty]$, 
as $n$ goes to infinity, are asymptotic to
$\sqrt{\pi/2}\, n^{1/2}$.
\end{corollary}

\begin{proof}
Apply the results of Theorem~\ref{th:mdx} and Theorem~\ref{th:mdy}
to the following inequalities
\begin{equation*}
D_n^y\le \max_{0\le i\le n}\Vert S_i\Vert_\infty 
\le \max_{0\le i\le n}\Vert S_i\Vert_1\le D_n^x+D^y_n.
\end{equation*}
\end{proof}

\subsection{Walk dimension}\label{ss:wd}

The maximal deviation of the walk in $n$ steps is linked to the
first exit time from a ball of radius $k$. Indeed if we put
$T_k=\min\{i:S_i\not\in B_k\}$, where $B_k$ is the ball of radius
$k$ centered in $(0,0)$, then
\begin{equation*}
\left( \max_{0\le i\le n}\Vert S_i\Vert\, \le k
\right) =\left(T_k\ge n\right).
\end{equation*}
Clearly the radius of the ball and $\Vert S_i\Vert$ must be computed
with respect to the same norm on the graph. We write $T^\infty_k$
and $T^1_k$ for the exit times with respect to the two norms we defined
in Section~\ref{sec:prelim}.

Recall now that given the simple random walk on a graph, if 
$\E[T_n]$ is of order $n^\alpha$, 
then by definition $\alpha$ is the walk dimension of the graph.
Usually, the norm with the respect to which the radius is computed
is $\Vert\cdot\Vert_1$, but as we will show, we may equivalently
consider $\Vert\cdot\Vert_\infty$.
Therefore we are now interested in the asymptotic behaviour of
$\E[T_n^\infty]$.

\begin{proposition}\label{th:tninfty}
$\E[T_n^\infty]\stackrel{n\to\infty}\sim n^2$.
\end{proposition}

\begin{proof}
We write
\begin{equation*}
\E[T_n^\infty]=\sum_{k\ge 0} \Prob(T_n^\infty>k),
\end{equation*}
that is, $\E[T_n^\infty]$ is equal to $\Theta_n(1)$,
where $\Theta_n(z)$
is the generating function of the sequence $\{\Prob\left(
T_n^\infty>k\right)\}_{k\ge0}$.
Let us observe that $\Prob(T_n^\infty>k)=\Prob(\max_{0\le i\le k}
\Vert S_i\Vert_\infty\le n)$.
We claim that
\begin{equation}\label{eq:Theta}
\Theta_n(z)  
=\frac{(1+w_n^2)(1-w_n^{n+1})^2}{(1-w_n)^2(1+w_n^{2n+2})}
\cdot\frac{(1+v^2)(1-v^{n+2})(1-v^{n+1})}{(1-v)(1-v^{2n+4})}
\end{equation}
where $v$ is such that \eqref{eq:v&z2} holds
and $w_n/(1+w_n^2)=v(1-v^{2n+2})/(2(1-v^{2n+4}))$.
Then since for $v=1$ (choosing the solution $w_n$ which is bounded 
in some neighbourhood of $v=0$)
\begin{equation*}
w_n=1+\frac{1-\sqrt{2n+3}}{n+1},
\end{equation*}
we get that 
\begin{equation*}
\Theta_n(1)\stackrel{n\to\infty}\sim n^2.
\end{equation*}
We are left with the proof of equation~\eqref{eq:Theta}.
We proceed as in Lemma~\ref{th:pdnx}, that is we separately consider
the two components of the walk.
Let $\widetilde G_n(0,l|z)$ be the generating function
of the probabilites of the 
$m$-step excursions along one single tooth, ending at a height $l$
or $-l$, with maximal deviation bounded by $n$ (see Figure~\ref{fig:absvert}). 
Put $\widetilde G_n(z):=\widetilde G_n(0,0|z)$. Then
\begin{equation*}
\Theta_n(z)=\psi_n\left(\frac{z\widetilde G_n(z)}4\right)
\sum_{l=0}^n\widetilde G_n(0,l|z),
\end{equation*}
where $\psi_n(z)$ is the generating function we already used in
Proposition~\ref{th:genfundxv}.
By \cite[Lemma 1.13]{cf:Woess}
\begin{equation*}
\widetilde G_n(z)=
\frac 1{1-z\widetilde F_n(1,0|z)/2},
\end{equation*}
where $\widetilde F_n(l,0|z)$ are the generating functions
of the $m$-step ``purely vertical'' excursions, from $(x,l)$
to $(x,0)$, with maximal deviation bounded by $n$.
The $\widetilde F_n$'s are the solutions of the linear system:
\begin{equation}
\mathbf A_{n-1}\cdot
\left[\begin{array}{c}
\widetilde F_{n}(1,0|z) \\
\widetilde F_{n}(2,0|z) \\
\vdots \\ \widetilde F_{n}(n,0|z) \\
\end{array}\right]
=
\left[\begin{array}{c}
z/2 \\ 0 \\
\vdots \\ 
0 \\
\end{array}\right].
\label{eq:ls}
\end{equation}
Then, with the substitution \eqref{eq:v&z2} (compare with \cite{cf:PP}), 
we get
\begin{equation*}
\begin{split}
\widetilde F_n(l,0|z)& =\frac{v^l(1-v^{2n-2l+2})}{(1-v^{2n+2})},\qquad l=1,\ldots,n,\\
\widetilde G_n(z)& =\frac{(1+v^2)(1-v^{2n+2})}{(1-v^{2n+4})}.
\end{split}
\end{equation*}
Thus, $z\widetilde G_n(z)/4=v(1-v^{2n+2})/(2(1-v^{2n+4}))=w_n/(1+w_n^2)$.
For the final excursion we must compute $\widetilde G_n(0,l|z)$, for
$l=1,\ldots n$. This can be done writing and solving 
$n$ linear systems in the spirit of
\eqref{eq:ls} to obtain
\begin{equation*}\begin{split}
\widetilde G_n(0,l|z) & =
\frac{v^l(1+v^2)(1-v^{2n-2l+2})}{(1-v^{2n+4})},\\
\sum_{l=0}^n \widetilde G_n(0,l|z) & =
\frac{(1+v^2)(1-v^{n+1})(1-v^{n+2})}{(1-v)(1-v^{2n+4})}.
\end{split}
\end{equation*}
\end{proof}

\begin{corollary}
$$\E[T_n^1]\stackrel{n\to\infty}\sim n^2.$$
\end{corollary}

\begin{proof}
Define $T_{n,\alpha}$ as the first exit time
from the rectangle $\{(x,y):|x|\le \alpha h,|y|\le (1-\alpha)h\}$,
and note that $T_{n,\alpha}\le T_n^1 \le T_n^\infty$.
A calculation similar to that of Proposition~\ref{th:tninfty}, 
shows that $\E[T_{n,\alpha}]\stackrel{n\to\infty}
\sim (1-\alpha)^2n^2$. Since $\alpha$ is an arbitrary
positive constant, we are done.
\end{proof}

\begin{corollary}
The walk dimension of 
$\mathbf C_2$ is equal to $2$.
\end{corollary}

%%%%%%%%%%%%%%%%%%%%%%%%%%%%%%%%%%%%%%%%%
%%%	SPAN        		      %%%
%%%%%%%%%%%%%%%%%%%%%%%%%%%%%%%%%%%%%%%%%
\section{Mean maximal span}
\label{sec:maxspan}

\begin{theorem}\label{th:msx}
$\E[M^x_n]\stackrel{n\to\infty}\sim \frac{2^{1/4}}{\Gamma(5/4)}n^{1/4}$.
\end{theorem}

\begin{proof}
Let $m^x_n=\max\{S_i^x\,:\,0\le i\le n\}$.
Then it is clear that $\E[M_n^x]=2\E[m^x_n]$.
Our first aim is to compute the generating function of
$\{\Prob(m^x_n\le h)\}_{n\ge0}$, which we denote by
$\widetilde \Psi_h(z)$. Then
\begin{equation*}
\widetilde \Psi_h(z)=\lim_{k\to\infty}\sum_{l=-k}^h
\widetilde\Psi_{h,k;l}(z)\cdot E(z),
\end{equation*}
where $\widetilde\Psi_{h,k;l}(z)$ is the generating function of the
probabilities of the $n$-step paths such that $-k\le S_i^x\le h$
for $0\le i\le n$, $S_n=(l,0)$ and $S_{n-1}^x\neq l$,
while $E(z)$ was defined and computed in Lemma~\ref{th:pdnx}.
Let us note that
$\widetilde\Psi_{h,k;l}(z)=\Psi_{h,k;l}(z\widetilde G(0,0|z)/4)$,
where $\Psi_{h,k;l}(z)$ is the generating function of the number
of the $n$-step paths on $\Z$ which stay in the interval between
$-k$ and $h$, and end at $l$. 
The functions $\Psi_{h,k;l}(w)$ are determined by
the linear system used in
\cite[Theorem 4.1]{cf:PP} to determine
$\Psi_{h,k;0}(w)$. Thus
\begin{equation*}
\begin{split}
\Psi_{h,k;l}(w) & =
\frac{w^l a_{h-l-1}(w)a_{k-1}(w)}{a_{h+k}(w)},\text{ if $l\ge0$},\\
\Psi_{h,k;l}(w) & =
\frac{w^{-l} a_{h-1}(w)a_{k+l-1}(w)}{a_{h+k}(w)},\text{ if $l\le -1$}.
\end{split}
\end{equation*}
Then we put $w=z\widetilde G(0,0|z)/4=v/(1+v^2)$ and we
obtain that (compare with the proof of Proposition~\ref{th:genfundxv})
\begin{equation*}\begin{split}
\lim_{k\to\infty}\sum_{l=-k}^h
\widetilde\Psi_{h,k;l}(z) & =\frac{(1+v^2)(1-v^{h+1})}{(1-v)^2},\\
E(z)&=\frac{1+6v^2+v^4}{(1+v^2)(1-v)^2}.
\end{split}
\end{equation*}
Then
\begin{equation*} \begin{split}
\widetilde\Psi_h(z)% & = \frac {1+6v^2+v^4}{(1-v)^4(1+v)}\left((1-v^{h+1})(1-v^{h+2})+v(1-v^{2h+2})\right)\\
&=\frac{(1+6v^2+v^4)(1-v^{h+1})}{(1-v)^4},\\
\sum_{n=0}^\infty \E[m_n^x]z^n&=
\frac{1+6v^2+v^4}{(1-v)^4}\sum_{h\ge0}v^{h+1}
\stackrel{z\ra1^-}\sim\frac{1}{2^{3/4}(1-z)^{5/4}},
\end{split}
\end{equation*}
whence $\E[m^x_n]\stackrel{n\to\infty}\sim 2^{-3/4}n^{1/4}/\Gamma(5/4)$
and we are done.
\end{proof}

\begin{theorem}\label{th:msy}
$\E[M^y_n]\stackrel{n\to\infty}\sim \sqrt{\frac8\pi}\,n^{1/2}.
$
\end{theorem}

\begin{proof}
Let $m^y_n=\max\{S_i^y:0\le i\le n\}$,
then $\E[M^y_n]=2\E[m^y_n]$. 
Denote by $\widehat\Psi_h(z)$ the generating function
of $\{\Prob\left(m_n^y\le h\right)\}_{n\ge0}$.
Then
\begin{equation*}
\widehat \Psi_h(z)=\lim_{k\to\infty}\sum_{l=-k}^h
\widehat \Psi_{h,k;l}(z),
\end{equation*}
where $\widehat \Psi_{h,k;l}(z)$ is the generating function of the
probabilities of the $n$-step paths (see Figure~\ref{fig:vertcomp}),
such that $-k\le S_i^y\le h$
for $0\le i\le n$, $S_n^y=l$. 
The functions $\widehat \Psi_{h,k;l}(z)$ are determined by the linear
system (where only the non-zero terms are displayed)
\begin{equation*}
\left[\begin{array}{rcl}
\left[\begin{array}{lcr}
\       & \                & \  \\
\       & \mathbf A_{h-1}  & \  \\
\       & \                & \  \\
\end{array}\right]
%& \ & \ \\
& \begin{array}{c}
\ \\ \  \\
 -z/4 \\
 \end{array}
 & \\ % fine parte superiore
   -z/2  &  (1-z/2) &    -z/2 \\ % riga centrale
   & \begin{array}{c}
 -z/4 \\ \ \\ \ \\
 \end{array}  &\left[
\begin{array}{lcr}
\       & \                & \  \\
\       & \mathbf A_{k-1}  & \  \\
\       & \                & \  \\
\end{array}\right]

\end{array}\right]
\left[\begin{array}{c}
\widehat \Psi_{h,k;h}(z) \\ \vdots
\\
\widehat \Psi_{h,k;1}(z) \\
\widehat \Psi_{h,k;0}(z) \\
\widehat \Psi_{h,k;-1}(z) \\
\vdots \\ \widehat \Psi_{h,k;-k}(z)
\end{array}\right]
=
\left[\begin{array}{c}
0\\ \vdots \\ 0 \\
1 \\ 0 \\
\vdots \\ 0 \\
\end{array}\right].
\end{equation*}
Denote by $b_{h,k}$ the determinant of the previous matrix.
For $h\ge2$ elementary computation lead to
\begin{equation*}
\widehat\Psi_{h}(z) =\lim_{k\to\infty}
\frac 1{2b_{h,k}}
\left(
a_{k-1}\sum_{l=0}^h\left(\frac z2\right)^la_{h-l-1}+
a_{h-1}\sum_{l=0}^k\left(\frac z2\right)^la_{k-l-1}
\right)
\end{equation*}
The interest in $\widehat\Psi_{h}(z)$ originates in that
\begin{equation}\label{eq:emy}
\sum_{n\ge0} \E[m_n^y]z^n=\sum_{h\ge0}\left(
\frac1{1-z}-\widehat\Psi_h(z)\right).
\end{equation}
By Lemma~\ref{th:ai} it is not difficult to prove that for $h\ge2$
(and $v$ such that \eqref{eq:v&z2} holds)
\begin{equation*}
\frac 1{1-z}-\Psi_h(z)=\frac{1+v^2}{(1-v)^2}
\frac{(1+v)v^{h+1}}{\left(2-(1-v)v^{2h+2}\right)}.
\end{equation*}
Noting that $2-(1-v)v^{2h+2}\in [2v,2]$ we get that
as $z\to1^-$ (and $v\to1^-$)
\begin{equation*}
\sum_{n\ge0}\E[m^y_n]z^n\sim\frac 2{(1-v)^3}\sim
\frac 1{\sqrt 2(1-z)^{3/2}}.
\end{equation*}
Of course we have to prove that the terms in \eqref{eq:emy}
corresponding to $h=0,1$ are negligible, but this follows
from elementary computation.
By \cite[Corollary 2]{cf:FO}
we deduce that $\E[m^y_n]\stackrel{n\to\infty}\sim\sqrt{\frac 2\pi}\,n^{1/2}$
and we are done.
\end{proof}

%%%%%%%%%%%%%%%%%%%%%%%%%%%%%%%%%%%%%%%%%%%%%%%%
%%% SCALING
%%%%%%%%%%%%%%%%%%%%%%%%%%%%%%%%%%%%%%%%%%%%%%%%

\section{Scaling limits}
\label{sec:scaling}

In the preceding sections we have seen that the expected values 
of the distances (with various meanings of this word) reached in $n$
steps are of order $n^{1/4}$ for the horizontal direction and of 
order $n^{1/2}$ for the vertical direction.
These results lead us to a natural question: what is the asymptotic
behaviour of the process where the horizontal component of the position
after $n$ steps is divided by $n^{1/4}$ and the vertical component is
divided by $n^{1/2}$? Of course we have to make this question more
precise.

In order to study the scaling of the process we choose a suitable
realization for the sequence $\{S_n\}_{n\ge0}$: 
let $X=\{X_n\}_{n\ge0}$ be
a sequence of random variables 
representing a simple random walk on $\Z$, 
and $Y=\{Y_n\}_{n\ge0}$ be a sequence
representing the random walk on $\Z$ moving according to Figure 2.
Choose $X$ and $Y$ to be independent and let $X_0=Y_0\equiv0$ a.s..
Moreover, let $L_k$ be the number of loops performed 
by $Y$ up to time $k$, that is
\[
L_k=\sum_{i=0}^{k-1}\ident_{Y_i=0,Y_{i+1}=0}.
\]
Clearly, $S_n=(X_{L_n},Y_n)$ is a realization 
of the position of the simple random walker on $\mathbf C_2$ at time $n$.
We are now able to define, by linear interpolation of the
three discrete time processes $X,Y$ and $L$,
a continuous time process $(X_{L_{nt}},Y_{nt})$.

\begin{theorem}\label{th:scaling}
\begin{equation}\label{eq:scaling}
\left(\frac {X_{L_{nt}}}{\sqrt[4]n},\frac{Y_{nt}}{\sqrt n}\right)_{t\ge0}
\stackrel{Law}\longrightarrow
\left(W_{L^0_t(B)},B_t\right)_{t\ge0},
\end{equation}
where $W$ and $B$ are two independent Brownian motions and
$L^0_t(B)$ is the local time at 0 of $B$.\label{th:scallim}
\end{theorem}

The theorem will be a consequence of Proposition~\ref{th:YL}
and Proposition~\ref{th:XL}.
We introduce the following notion of convergence of 
stochastic processes (see Definition 2.2 of \cite{cf:CSY}).
\begin{definition}
A sequence of $\R^k$-valued
stochastic processes $(Z_t^n; t\ge0)_{n\ge0}$ converges to
a process $(Z_t; t\ge0)$ in probability uniformly on compact intervals
if for all $t\ge0$, as $n\to\infty$
\[
\sup_{s\le t}\Vert Z_s^n-Z_s\Vert\stackrel{\Prob}{
\longrightarrow}0,
\]
where $\Vert\cdot\Vert$ is a norm on $\R^k$ (for instance
$\Vert\cdot\Vert_1$).
We will briefly write
\[
(Z_t^n; t\ge0)\stackrel{U.P.}{
\longrightarrow}(Z_t; t\ge0).
\]
\end{definition}
Since U.P. convergence of a vector is equivalent to
U.P. convergence of its components and implies convergence in distribution,
in order to prove Theorem~\ref{th:scallim} it will suffice to prove
that each component in \eqref{eq:scaling} U.P. converges to the corresponding
limit.

The main idea is that $Y$ and $L$ are not much different from,
respectively, a simple random walk $Y^\prime$ on $\Z$ and the process
$L^\prime$ which counts the visits of $Y^\prime$ to 0.
There is a natural correspondence between $Y$ and $Y^\prime$, so let
us define $Y^\prime$ and some other auxiliary variables which will be
needed in the sequel.
Given $Y$, let $l$ and $N$ be the processes which respectively
count its returns to 0  (not including the loops and counting time 0 as the
first ``return'') and 
the time spent not looping at 0. Namely, let $l_0=1$, and
$l_k =1+\sum_{i=0}^{k-1} \ident_{Y_{i+1}=0,Y_{i}\neq0}$ 
 and $N_k=\sum_{i=0}^{k-1}\ident_{Y_i\neq Y_{i+1}}$
for $k\ge1$. Clearly $N_k=k-L_k$. Moreover, note that for $k\ge1$,
\begin{equation}\label{eq:Lk}
\sum_{i=0}^{l_k-1}\tau_i\le L_k\le\sum_{i=0}^{l_k}\tau_i,
\end{equation}
where $\tau=\{\tau_i\}_{i\ge0}$ is a suitable
sequence of iid random variables with
geometric distribution of parameter 1/2.
Now define a simple random walk $Y^\prime$ on $\Z$
by $Y_n=Y^\prime_{N_n}$. % ($Y^\prime$ does not take into account the time spent by $Y$ looping at 0). 
Then $L^\prime_k=\sum_{i=0}^{k}\ident_{Y_i^\prime =0}$ ($L^\prime$
counts the visits at 0 or, equivalently, the returns to 0).
We note that $l_n  = L^\prime_{N_n}$ and $0\le l_n\le L^\prime_n$.
We first prove a property of $L$.

\begin{lemma}\label{th:limLn}
\[\frac{L_n}{\sqrt n}\stackrel{Law}\lra |\mathcal N(0,1)|.\]
\end{lemma}
\begin{proof}
The main ideas of the proof are the facts that
$L^\prime_n/\sqrt n\to |\mathcal N(0,1)|$ and that $N_n$ is 
not much different from $n$.
Indeed one can easily prove the first fact
(for the distribution of $(L^\prime_{2n}-1)$ see 
\cite[Chapter III, Exercise 10]{cf:Feller}).
%
%The proof of the first fact is nothing more than an exercise:
%from \cite[Chapter III, Exercise 10]{cf:Feller} 
%\begin{equation*}
%\Prob\left(L^\prime_{2n}=r\right)=\frac 1{2^{2n-r}}\binom{2n-r}{n}.
%\end{equation*}
%Then, using Stirling's formula,
%\begin{equation*}
%\Prob\left(\frac {L^\prime_{2n}}{\sqrt{2n}}\le x\right)=
%\sum_{r:\frac r{\sqrt{2n}}\le x}
%\frac 1{2^{2n-r}}\binom{2n-r}{n}\stackrel{n\to\infty}\lra
%\sqrt{\frac 2\pi}\int_0^x
%\exp(-t^2/2)dt,
%\end{equation*}
%that is, the limit is the distribution function of $|\mathcal N(0,1)|$.
%The statement for odd times follows from $|L^\prime_{2n+1}-L^\prime_{2n}|\le1$.
%
%As for $N_n$, note that $n-\sum_{i=1}^{L^\prime_n}\tau_i\le N_n\le n$,
%and since $L^\prime$ and $\tau$ are independent,
%\begin{equation*}
%\frac{\sum_{i=0}^{L^\prime_n}\tau_i}n\stackrel{a.s.}\ra 0.
%\end{equation*}
%whence $N_n/n\stackrel{a.s.}\to1$. 
%
By \eqref{eq:Lk}, the thesis is a consequence of
\begin{equation}\label{eq:taui}
\frac{\sum_{i=0}^{l_n}\tau_i}{\sqrt n}
\stackrel{Law}{\to} |\mathcal N(0,1)|.
\end{equation}
By the strong law of large numbers and Slutzky's theorem,
\begin{equation*}
\frac{\sum_{i=0}^{L^\prime_n}\tau_i}{\sqrt n}
=\frac{\sum_{i=0}^{L^\prime_n}\tau_i}{L^\prime_n}
\cdot\frac{L^\prime_n}{\sqrt n}
\stackrel{Law}{\to} |\mathcal N(0,1)|.
\end{equation*}
Then %, since $\sum_{i=0}^{l_n}\tau_i\le\sum_{i=0}^{L^\prime_n}\tau_i$,
\eqref{eq:taui} will follow once we show that
$\sum_{i=0}^{l_n}\tau_i/ \sum_{i=0}^{L^\prime_n}\tau_i\stackrel\Prob\to1$. 
%where $\{V_n\}_{n\ge0}$ is a sequence of random variables such that $V_n/\sqrt n\stackrel{Law}{\to} |\mathcal N(0,1)|$.
Indeed note that
\begin{equation}\label{eq:ln}
L^\prime_n=L^\prime_{N_n}+R^\prime_{n} =l_n+R^\prime_n,
\end{equation}
where $R^\prime_n$ is the number of visits to 0 of $Y^\prime$
between time $N_n$ and time $n$.
Let $T_n$ be the time $Y^\prime$ first visits 0 after 
time $N_n$,
and for any $j\ge0$ let $L^{''}_j$ be the number of visits to 0 
before time $j$ of the random walk $Y^{''}_m:=
Y^\prime_{m+T_n}$. Clearly $L^{''}_j$ is independent
of $\{Y^\prime_k\}_{k=0}^{N_n}$ and has the same distribution
of $L^\prime_j$.
Then 
\[
\frac{\sum_{i=0}^{l_n}\tau_i}{\sum_{i=0}^{L^\prime_n}\tau_i}=
1-\frac{\sum_{i=l_n+1}^{L^\prime_n}\tau_i}{\sum_{i=0}^{L^\prime_n}\tau_i}
=
1-\frac{\sum_{i=1}^{R^\prime_n}\widetilde\tau_i}{\sum_{i=0}^{L^\prime_n}
\tau_i},
\]
where $\widetilde\tau_i=\tau_{i+l_n}$ and $\sum_{i=1}^{R^\prime_n}
\widetilde\tau_i$
is equal to zero if $L^\prime_n=l_n$. We are left with the proof that
${\sum_{i=1}^{R^\prime_n}\widetilde\tau_i}/{\sum_{i=0}^{L^\prime_n}
\tau_i}\stackrel\Prob\to0$. 
Since $0\le R^\prime_n\le L^{''}_{n-N_n}$, it suffices to prove that
${\sum_{i=1}^{L^{''}_{n-N_n}}\widetilde\tau_i}/{\sum_{i=0}^{L^\prime_n}
\tau_i}\stackrel\Prob\to0$.
Fix $\eps>0$: since $\sum_{i=0}^{L^\prime_n}\tau_i\ge n-N_n$ we have that
\begin{equation*}\begin{split}
\Prob\left(\frac{\sum_{i=1}^{L^{''}_{n-N_n}}
\widetilde\tau_i}{\sum_{i=0}^{L^\prime_n}\tau_i}>\eps
\right)&\le \sum_{k=1}^{n}\Prob\left(\frac{\sum_{i=1}^{L^{''}_{n-N_n}}
\widetilde\tau_i}{n-N_n}>\eps,\ n-N_n=k
\right)\\
%& \le \sum_{k=1}^{n}\Prob\left(\frac{\sum_{i=1}^{L^{''}_{k}}
%\widetilde\tau_i}{k}>\eps,\ n-N_n=k\right)\\
&=\sum_{k=1}^{n}\Prob\left(\frac{\sum_{i=1}^{L^{''}_{k}}
\widetilde\tau_i}{k}>\eps\right)\Prob(n-N_n=k).
\end{split}
\end{equation*}
Now fix $\delta>0$ and choose $M=M(\delta)$ such that 
$\Prob(\sum_{i=1}^{L^{''}_{k}}
\widetilde\tau_i/k>\eps)<\delta$ for all $k\ge M$ (this is possible
by the law of large numbers using the facts that
$L^{''}$ and $\widetilde\tau$ are independent and 
$L^{''}_k/k\stackrel\Prob\to0$). Then
\begin{equation*}
\sum_{k=1}^{n}\Prob\left(\frac{\sum_{i=1}^{L^{''}_{k}}
\widetilde\tau_i}{k}>\eps\right)\Prob(n-N_n=k)
\le \Prob(n-N_n<M)+\delta.
\end{equation*}
Since $n-N_n\stackrel{\Prob}\to\infty$ we are done.
%if we put $L^{''}_{j}=0$ if $j<0$,we have that\begin{equation}\label{eq:ln}L^\prime_n=L^\prime_{N_n}+L^{''}_{n-T_n} \le L^\prime_{N_n}+
%L^{''}_{n-N_n}=l_n+L^{''}_{n-N_n},\end{equation}
%hence \begin{equation*}\sum_{i=0}^{l_n}\tau_i\ge\sum_{i=0}^{L^\prime_n(1-L^{''}_{n-N_n}/L^\prime_n)}\tau_i,\end{equation*}and \eqref{eq:taui} follows once we show that $L^{''}_{n-N_n}/L^\prime_n\stackrel{\Prob}\to0$. Write $\Prob\left(L^{''}_{n-N_n}\ge\delta L_n^\prime\right)$ as \begin{equation*}\begin{split} \sum_{j=0}^{n}& \Prob\left( L^{''}_{n-j}\ge\delta L^\prime_{n}\right) \Prob \left(N_n=j\right)\\&\le \Prob \left(N_n\le n(1-\eps)\right)+\sum_{j=\lfloor n(1-\eps)\rfloor}^{n}\Prob\left( L^{''}_{n-j}\ge\delta L^\prime_{n}\right) \Prob \left(N_n=j\right)\\&\le \Prob \left(N_n\le n(1-\eps)\right)+\Prob\left( L^{''}_{n-\lfloor n(1-\eps)\rfloor}\ge\delta L^\prime_{n}\right) .\end{split}\end{equation*}
%The first term, given $\eps>0$, is arbitrarily small if $n$ is sufficiently large ($N_n/n\stackrel{a.s.}0$),while the second term may be written as\begin{equation}\label{eq:pbde}
%\Prob\left(\left(\frac1{\sqrt n}L^\prime_n,\frac1{\sqrt{n\eps}}{L^{''}_{n-\lfloor n(1-\eps)\rfloor}}\right)\in B_{\delta,\eps}\right)\stackrel{n\to\infty}\lra\Prob\left( (|Z^\prime|,|Z^{''}|)\in B_{\delta,\eps}\right),\end{equation}where $B_{\delta,\eps}=\{(x,y)\in \R^+\times \R^+\, :\, y\ge x\delta/\sqrt{\eps}\}$ and $Z^\prime$ and $Z^{''}$are two independent standard normals. It is clearthat the probability in \eqref{eq:pbde} tends to zero as$\eps\to0^+$, and we are done.

\end{proof}
\begin{proposition}\label{th:YL}
\[
\left(\frac1{\sqrt n}{Y_{nt}},\frac1{\sqrt n}L_{nt}\right)_{t\ge0}
\stackrel{U.P.}\longrightarrow
\left(B_t,{L^0_t(B)}\right)_{t\ge0}.
\]
\end{proposition}
\begin{proof}

Consider the processes $Y^\prime$ and $L^\prime$ defined
before Lemma \ref{th:limLn} and  by interpolation define
the sequence of two-dimensional continuous time processes 
$\left((\frac1{\sqrt n}Y^\prime_{nt},\frac1{\sqrt n}L^\prime_{nt}) ;
t\ge0\right)_{n\ge0}$.
Then by Theorem 3.1 of \cite{cf:CSY} we have that
\[
\left(\frac1{\sqrt n}{Y^\prime_{nt}},\frac1{\sqrt n}L^\prime_{nt}\right)_{t\ge0}
\stackrel{U.P.}\longrightarrow
\left(B_t,{L^0_t(B)}\right)_{t\ge0}.
\]
To prove our statement, it suffices to show that these two properties
hold:
\[\begin{split}
(A): &\left(\frac1{\sqrt n}\left(Y_{\lfloor nt\rfloor}-
Y^\prime_{\lfloor nt\rfloor}\right)\right)_{t\ge0}
\stackrel{U.P.}\longrightarrow0\\
(B): &\left(\frac1{\sqrt n}\left(L_{\lfloor nt\rfloor}-
L^\prime_{\lfloor nt\rfloor}\right)\right)_{t\ge0}
\stackrel{U.P.}\longrightarrow0.
\end{split}
\]
Note that $Y^\prime$ is the sum of iid increments $\xi_i$ such that
$\Prob(\xi_i=\pm1)=1/2$, hence
\[\begin{split}
Y_\nt& =Y^\prime_{\nt-L_\nt}=
\sum_{i=1}^{\nt-L_\nt}\xi_i\\
&= \sum_{i=1}^{\nt}\xi_i-\sum_{i=\nt-L_\nt+1}^{\nt}\xi_i
=Y^\prime_\nt-\sum_{i=1}^{L_\nt}\tilde\xi_i,
\end{split}
\]
where $\tilde\xi_i=\xi_{\nt-L_\nt+i}$ (and $\sum_{i=1}^{L_\nt}\tilde\xi_i=0$
if $L_\nt=0$).
Thus $Y_{\lfloor nt\rfloor}-Y^\prime_{\lfloor nt\rfloor}=
\sum_{i=1}^{L_\nt}\tilde\xi_i$.
Then $(A)$ follows from  
\begin{equation}\label{eq:xitilde}
\sup_{s\le t}\left|\frac1{\sqrt n}\sum_{i=1}^{L_\nns}\tilde\xi_i
\right|\stackrel\Prob\longrightarrow0.
\end{equation}
Indeed
\[
\sup_{s\le t}\left|\sum_{i=1}^{L_\nns}\tilde\xi_i
\right|\le
\max_{k\le L_\nt}\left|\sum_{i=1}^{k}\tilde\xi_i
\right|,
\]
and if we denote by $M_n=\max_{k\le n}\sum_{i=1}^{k}\tilde\xi_i$
and by $m_n=\min_{k\le n}\sum_{i=1}^{k}\tilde\xi_i$, clearly 
$M_n$ and $-m_n$ are identically distributed, and
$\max_{k\le L_\nt}\left|\sum_{i=1}^{k}\tilde\xi_i
\right|=\max\left\{M_{L_\nt},-m_{L_\nt}\right\}$ 
Hence to prove \eqref{eq:xitilde} it suffices to show that
\[
\frac1{\sqrt n}M_{L_\nt}\stackrel\Prob\longrightarrow0.
\]
The distribution of $M_n$ is well known (see \cite[Chapter III.7]{cf:Feller}),
and it is easy
to show that $M_n/\sqrt n\stackrel{Law}\to |\mathcal N(0,1)|$.
Noting that $L_\nt$ is independent of $M_k$, we have that
\[\begin{split}
\Prob(M_{L_\nt}>\eps\sqrt n)&=
\sum_{k=0}^\nt
\Prob(M_{k}>\eps\sqrt n)\Prob(L_\nt=k)\\
&\le 
\sum_{\eps\sqrt n<k\le\alpha\sqrt\nt}
\Prob(M_{k}>\eps\sqrt n)\Prob(L_\nt=k)+
\sum_{k>\alpha\sqrt\nt}\Prob(L_\nt=k)\\
&\le \Prob(M_{\lfloor\alpha\sqrt\nt\rfloor}>\eps\sqrt n)+
\Prob(L_\nt>\alpha\sqrt\nt).
\end{split}
\]
By Lemma \ref{th:limLn}, for any positive $\eps^\prime$ and $t$ there exist
$\alpha$ and $n^\prime$ such that
$\Prob(L_\nt>\alpha\sqrt\nt)<\eps^\prime$ for all $n\ge n^\prime$.
On the other hand for any given $\eps,\alpha$
\[
\Prob\left(\frac{M_{\lfloor\alpha\sqrt\nt\rfloor}}{\sqrt\alpha\sqrt[4]\nt}>
\frac{\eps\sqrt[4]n}{\sqrt\alpha\sqrt[4] t}\right)
\stackrel{n\to\infty}\lra 0,
\]
whence $(A)$ is proven.
\medskip

%%%%%
Now, let us address to $(B)$. We first note that we may consider
a mapping between the number of steps taken by $Y^\prime$ and the ones 
taken by $Y$. Indeed when $Y^\prime$ has taken $\nt$ steps, 
then $Y$ has taken $\nt+\sum_{i\le L^\prime_\nt}
\tau_i$ steps
(that is, if $Y^\prime_\nt=0$ we decide to count for $Y$ all the loops
it performs after this last return to 0).
Let us write
\begin{multline}
\frac1{\sqrt n}\left(
L^\prime_\nt-L_\nt\right)=
\frac1{\sqrt n}\left(
L^\prime_\nt-L_{\nt+\sum_{i\le L^\prime_\nt}\tau_i}\right)+\\
\frac1{\sqrt n}\left(
L_{\nt+\sum_{i\le L^\prime_\nt}\tau_i}-L_\nt\right)=:I+II.
\end{multline}
We prove that $I\stackrel{U.P.}\to0$. Indeed
$L_{\nt+\sum_{i\le L^\prime_\nt}\tau_i}=\sum_{i\le L^\prime_\nt}\tau_i$, thus
\[
\sup_{s\le t}\left|
L^\prime_\nns-\sum_{i\le L^\prime_\nns}\tau_i
\right|=\sup_{s\le t}\left|\sum_{i\le L^\prime_\nns}(1-\tau_i)
\right|\le\max_{k\le L^\prime_\nt }\left|\sum_{i\le k}(1-\tau_i)
\right|.
\]
Now choose $\delta>0$. By independence of $L^\prime$ and $\tau$ we have that
the probability that \break
$\max_{k\le L^\prime_\nt }\left|\sum_{i\le k}(1-\tau_i)
\right|$
is larger than $\eps\sqrt n$ is bounded by
\[%begin{multline*}
%\Prob\left(\sup_{0\le k\le L^\prime_\nt }\left|\sum_{i\le k}(1-\tau_i)
%\right|>\eps\sqrt n\right) =
%\sum_{j=1}^{\nt /2}\Prob\left(\sup_{0\le k\le j}\left|\sum_{i\le k}(1-\tau_i)
%\right|>\eps\sqrt n\right)\Prob\left(L^\prime_\nt=j\right)\\
%\le
\Prob\left(L^\prime_\nt>\alpha\sqrt\nt\right)+
\sum_{j\le\alpha\sqrt\nt}\Prob\left(\max_{k\le j}\left|
\sum_{i\le k}(1-\tau_i)
\right|>\eps\sqrt n\right)\Prob\left(L^\prime_\nt=j\right).
\]%\end{multline*}
The first term is smaller than $\delta$ if
%$\alpha\ge\alpha(\delta)$ and $n\ge n(\alpha,\delta)$.
$\alpha$ and $n$ are sufficiently large.
As for the second term, it is clearly less or equal to 
\begin{equation}\label{eq:probtau1}
\Prob\left(\max_{k\le \alpha\sqrt\nt}\left|
\sum_{i\le k}(1-\tau_i)
\right|>\eps\sqrt n\right).
\end{equation}
Observe that, by the law of large numbers,
for any positive $\eps^\prime$ and $\delta$
there exists $k_0=k_0(\eps^\prime,\delta)$ such that for all $k\ge k_0$
\[
\Prob\left(
\left|
\frac{\sum_{i\le k}(1-\tau_i)}k
\right|<\eps^\prime
\right)\ge 1-\delta.
\]
Hence \eqref{eq:probtau1} is less or equal to
\[
\Prob\left(\max_{k\le k_0}\left|
\sum_{i\le k}(1-\tau_i)
\right|>\eps\sqrt n\right)+
\Prob\left(\max_{k_0\le k\le \alpha\sqrt\nt}\left|
\sum_{i\le k}(1-\tau_i)
\right|>\eps\sqrt n\right).
\]
The first term clearly tends to 0 as $n$ grows to infinity, while
the second term is not larger than
\[
\delta+\Prob\left(\max_{k_0\le k\le \alpha\sqrt\nt} \left|
{\sum_{i\le k}(1-\tau_i)} 
\right|>\eps\sqrt n,
\left|
\frac{\sum_{i\le k}(1-\tau_i)}k
\right|<\eps^\prime\quad\forall k\ge k_0
\right).
\]
But if $|\sum_{i\le k}(1-\tau_i)|/k<\eps^\prime$ for all $k\ge k_0$, then
\[
\sup_{k_0\le k\le \alpha\sqrt\nt} \left|
{\sum_{i\le k}(1-\tau_i)} 
\right|<\alpha\eps^\prime\sqrt\nt,
\]
which is smaller than $\eps\sqrt n$ if $\eps^\prime$ is sufficiently small.
This proves that  $I\stackrel{U.P.}\to0$.
%%%%%%%%%%%%%%%%%%%%%

\smallskip
We now prove that $II\stackrel{U.P.}\to0$. Indeed
\[
0\le L_{\nns+\sum_{i\le L^\prime_\nns}\tau_i}-L_\nns
\le\sum_{i\le L^\prime_\nns}\tau_i-\sum_{i\le l_\nns -1}\tau_i
=\sum_{i=l_\nns}^{L^\prime_\nns}\tau_i.\]
Using \eqref{eq:ln} and the definitions of $L^{''}$ and $\widetilde\tau$
thereafter, we have that
$\sum_{i=l_\nns}^{L^\prime_\nns}\tau_i\le
\sum_{i=0}^{L^{''}_{\nns-N_\nns}}\widetilde\tau_i$, and
% $l_\nns\ge L^\prime_\nns-L^{''}_{\nns-N_\nns}$, putting 
% $\tilde\tau_i=\tau_j$ with $j=i+L^\prime_\nns-L^{''}_{\nns-N_\nns}$,
for any positive $\eps^\prime$,
\[%\begin{split}
\sup_{s\le t}\sum_{i=0}^{L^{''}_{\nns-N_\nns}}\tilde\tau_i
=\max\left\{
\sup_{s\le \eps^\prime}\sum_{i=0}^{L^{''}_{\nns-N_\nns}}\tilde\tau_i,
\sup_{\eps^\prime<s\le t}\sum_{i=0}^{L^{''}_{\nns-N_\nns}}\tilde\tau_i
\right\}=:\max (C,D).%\end{split}
\]
Choose $\delta>0$.
Let us show that $\Prob(C>\eps\sqrt n)<\delta$  if
$\eps^\prime$ is sufficiently small and $n$ sufficiently large.
Indeed
\[
\sup_{s\le \eps^\prime}\sum_{i=0}^{L^{''}_{\nns-N_\nns}}\tilde\tau_i
\le \sum_{i=0}^{L^{''}_{\lfloor n\eps^\prime\rfloor}}\tilde\tau_i,
\]
and by independence of $L^{''}_{\lfloor n\eps^\prime\rfloor}$ and
$\tilde\tau$, %for $\alpha\ge\alpha(\delta)$, $\eps^\prime\le \eps^\prime
%(\alpha,\delta,\eps)$ and $n\ge n(\eps,\eps^\prime,\delta,\alpha)$,
\begin{equation}\label{eq:estc}
%\begin{split}
\Prob\left(\sum_{i\le L^{''}_{\lfloor n\eps^\prime\rfloor}}
\tilde\tau_i>\eps\sqrt n\right)
%&
\le \Prob\left(L^{''}_{\lfloor n\eps^\prime\rfloor}>\alpha
\sqrt{\lfloor n\eps^\prime\rfloor}\right)+
\Prob\left(\sum_{i\le\alpha\sqrt{\lfloor n\eps^\prime\rfloor}}
\widetilde\tau_i>\eps\sqrt n\right).
%\\&\le\delta+
%\Prob\left(\frac{\sum_{i\le\alpha\sqrt{\lfloor n\eps^\prime\rfloor}}
%\widetilde\tau_i}{\alpha\sqrt{n\eps^\prime}}>
%\frac\eps{\alpha\sqrt{\eps^\prime}}\right).%<2\delta.
%\end{split}
\end{equation}
Choose $\alpha$ such that 
$\Prob\left(L^{''}_{\lfloor n\eps^\prime\rfloor}>\alpha
\sqrt{\lfloor n\eps^\prime\rfloor}\right)<\delta$ for $n$
sufficiently large, and then $\eps^\prime$ sufficiently small
such that $\Prob\left(\sum_{i\le\alpha\sqrt{\lfloor n\eps^\prime\rfloor}}
\widetilde\tau_i>\eps\sqrt n\right)<\delta$.

\smallskip

Now keep $\eps^\prime$ fixed. In order to prove that
$\Prob(D>\eps\sqrt n)<3\delta$ for $n$ sufficiently large,
we observe that for any positive $\eps^{''}$
there exists $\bar n$ such that for all $n\ge \bar n$,
$\Prob\left(N_n\ge n(1-\eps^{''})\right)$ is larger than $1-\delta$.
Then for all $n$ such that $\lfloor n\eps^\prime\rfloor\ge\bar n$,
$\sup_{\eps^\prime< s\le t}\left(L^{''}_{\nns-N_\nns}\right)
\le L^{''}_{\lfloor nt\eps^{''}\rfloor}$
with probability larger than $1-\delta$. Hence, as in \eqref{eq:estc}
we have that, for $\alpha$ and $n$ sufficiently large and
$\eps^{''}$ sufficiently small.
\[\begin{split}
\Prob\left(D>\eps\sqrt n\right)
&\le\delta+
\Prob\left(
\sum_{i\le L^{''}_{\lfloor nt\eps^{''}\rfloor}}\tilde\tau_i>\eps\sqrt n\right)\\
&\le\delta+\Prob\left(L^{''}_{\lfloor nt\eps^{''}\rfloor}>\alpha
\sqrt{\lfloor nt\eps^{''}\rfloor}\right)+\Prob\left(
\sum_{i\le\alpha \sqrt{\lfloor nt\eps^{''}\rfloor}}
\tilde\tau_i>\eps\sqrt n\right)
<3\delta.\\
%&\le2\delta+\Prob(\frac{\sum_{i=1}^{\alpha\sqrt{\lfloor nt\eps^{''}\rfloor}}\tau_i}{\alpha\sqrt{nt\eps^{''}}}
%>\frac\eps{\alpha\sqrt{t\eps^{''}}})<3\delta,
\end{split}\]

%% is equal to the number of loops that $Y$ performs between time $\nns$ and time $\nns+\sum_{i\le L^\prime_\nns}\tau_i$
% which is not larger than $L_{\sum_{i\le L^\prime_\nns}\tau_i}$. Hence, since $L^\prime_\nns$ is an increasing function of $s$
% \[\sup_{0\le s\le t}\left|L_{\nns+\sum_{i\le L^\prime_\nns}\tau_i}-L_\nns\right|\le \sup_{0\le s\le t}L_{\sum_{i\le L^\prime_\nns}\tau_i}=L_{\sum_{i\le L^\prime_\nt}\tau_i}.\]
% Thus\[\begin{split}\Prob\left(L_{\sum_{i\le L^\prime_\nt}\tau_i}>\eps\sqrt n\right)&\le\sum_{\eps\sqrt n\le k\le\alpha\sqrt\nt}\Prob(L^\prime_\nt=k)\Prob(L_{\sum_{i\le k}\tau_i}>\eps\sqrt n)\\&+\Prob(L^\prime_\nt>\alpha\sqrt\nt)\\&\le\Prob(L_{\sum_{i\le \alpha\sqrt\nt}\tau_i}>\eps\sqrt n)+\Prob(L^\prime_\nt>\alpha\sqrt\nt).\end{split}\]The second term can be taken arbitrarily small if $\alpha$ is sufficiently large. As for the first term, it may be written as%\[\Prob\left(\frac{L_{\sum_{i\le \alpha\sqrt\nt}\tau_i}}{\sqrt{\sum_{i\le \alpha\sqrt\nt}\tau_i}}\sqrt{\frac{\sum_{i\le \alpha\sqrt\nt}\tau_i}{\alpha\sqrt{nt}}}\sqrt{\frac{\alpha\sqrt t}{\sqrt n}}>\eps\right).\]The first factor converges, a.s. with respect to $\{\tau_i\}_i$,to $|\mathcal N(0,1)|$, while the rest of the product converges a.s. to 0. % VEDERE!!!

\end{proof}

\begin{proposition}\label{th:XL}
\[
\left(\frac {X_{L_{nt}}}{\sqrt[4]n}\right)_{t\ge0}
\stackrel{U.P.}\longrightarrow
\left(W_{L^0_t(B)}\right)_{t\ge0}.
\]
\end{proposition}
\begin{proof}
Clearly the statement may be rephrased as
\[
\left(\frac {X_{L_{n^2t}}}{\sqrt n}\right)_{t\ge0}
\stackrel{U.P.}\longrightarrow
\left(W_{L^0_t(B)}\right)_{t\ge0}.
\]
Note that (writing $L^0_s$ instead of $L^0_s(B)$)
\[
\sup_{s\le t}\left|
\frac{X_{L_{n^2s}}}{\sqrt n}-W_{L^0_s}
\right|\le
\sup_{s\le t}\left|
\frac{X_{L_{n^2s}}-X_{nL^0_s}}{\sqrt n}
\right|+
\sup_{s\le t}\left|
\frac{X_{nL^0_s}}{\sqrt n}-W_{L^0_s}
\right|=(A)+(B).
\]
We have to show that $(A)$ and $(B)$ $\stackrel{U.P.}\to0$: it
suffices to prove the statement for $L_{\lfloor n^2s\rfloor}$
and $\lfloor nL^0_s\rfloor$ instead of $L_{n^2s}$
and $nL^0_s$ respectively.
Represent $X_k=\sum_1^k\widetilde\eta_i$, where $\widetilde\eta
=\{\widetilde\eta_i\}_{i\ge1}$ is 
an iid family such that
$\Prob(\widetilde\eta_i=\pm1)=1/2$, then
\[
\sup_{s\le t}\left|
{X_{L_{\lfloor n^2s\rfloor}}-X_{\lfloor nL^0_s\rfloor}}
\right|=\sup_{s\le t}\left|
{\sum_{i=1}^{n\Delta_{n,s}}\eta_i}
\right|,
\]
where $n\Delta_{n,s}:=|L_{\lfloor n^2s\rfloor}-\lfloor nL^0_s\rfloor|$
and $\eta_i=\widetilde\eta_{\min\{L_{\lfloor n^2s\rfloor},
\lfloor nL^0_s\rfloor\}+i}$.
Note that $(\Delta_{n,t})_{t\ge0}\stackrel{U.P.}\to0$, indeed
\[
\Delta_{n,t}\le
\left|\frac {L_{\lfloor n^2t\rfloor}}n-L^0_t\right|
+\left|L^0_t-\frac{\lfloor nL^0_t\rfloor}n\right|
\le\left|\frac {L_{\lfloor n^2t\rfloor}}n-L^0_t\right|+\frac1n,
\]
and both these summands tends U.P. to 0 (the first by Proposition \ref{th:YL}).
In order to show that
\[
\Prob\left(
\sup_{s\le t}
\left|\sum_{i=1}^{n\Delta_{n,s}}\eta_i\right|>\eps\sqrt n
\right)\stackrel{n\to\infty}\to0
\]
it suffices to prove that
\[
\Prob\left(
\max_{k\le\sup_{s\le t}(n\Delta_{n,s})}
\sum_{i=1}^{k}\eta_i>\eps\sqrt n
\right)\stackrel{n\to\infty}\to0.
\]
Now put $M_n=\max_{k\le n}\sum_{i=1}^{k}\eta_i$ and recall that
$(n\Delta_{n,s})_{n\ge0,s\le t}$ and $\eta$ are independent:
this last probability may be written as
\[\begin{split}
&\sum_{j=0}^\infty\Prob(M_j>\eps\sqrt n)\Prob\left(\sup_{s\le t}
(n\Delta_{n,s})=j\right)\\
%&\le \Prob(\sup_{s\le t}(\Delta_{n,s})>\alpha)+\sum_{j\le\alpha n}\Prob(M_j>\eps\sqrt n)\Prob(\sup_{s\le t}(n\Delta_{n,s})=j)\\
&
\le \Prob\left(\sup_{s\le t}
(\Delta_{n,s})>\alpha\right)
+ \Prob(M_{\lfloor n\alpha\rfloor}>\eps\sqrt n).
\end{split}\]
For all positive $\alpha$, 
$\Prob\left(\sup_{s\le t}
(\Delta_{n,s})>\alpha\right)$
can be made arbitrarily
small if $n$ is picked large enough, while 
\[
\Prob\left({M_{\lfloor n\alpha\rfloor}}>\eps\sqrt n\right)
\stackrel{n\to\infty}\to\Prob\left(|\mathcal N(0,1)|>\frac\eps{\sqrt\alpha}
\right).
\]
Thus if $\delta>0$ is fixed, by choosing
$\alpha$ sufficiently small we get
$\Prob\left({M_{\lfloor n\alpha\rfloor}}>\eps\sqrt n\right)<\delta$
if $n$ is sufficiently
large. This proves $(A)\stackrel{U.P.}\to0$.
\smallskip

Now note that $X_n$ and $L^0_s$ are independent, so we may think of
$X_{\lfloor nL^0_s\rfloor}$ 
as defined on a product probability space $\Omega_1\otimes
\Omega_2$ (that is $X_{\lfloor nL^0_s\rfloor}(\omega_1,\omega_2)=
X_{\lfloor nL^0_s(\omega_2)\rfloor}(\omega_1)$). 
For any fixed $\omega_2$ we have that 
\begin{equation}\label{eq:probx1}
\Prob_1\left(
\sup_{s\le t}\left|\frac{X_{\lfloor nL^0_s(\omega_2)\rfloor}}{\sqrt n}-
W_{L^0_s(\omega_2)}\right|>\eps
\right)\stackrel{n\to\infty}\to0.
\end{equation}
Indeed it is known (see \cite{cf:CSY}) that 
$\Prob_1\left(
\sup_{s\le t}\left|\frac{X_{ns}}{\sqrt n}-
W_{s}\right|>\eps
\right)\stackrel{n\to\infty}\to0$,
and since for any fixed $\omega_2$, as $s\uparrow t$,
$L^0_s(\omega_2)\uparrow L_t^0(\omega_2)$,
\[
\sup_{s\le t}\left|\frac{X_{\lfloor nL^0_s(\omega_2)\rfloor}}{\sqrt n}-
W_{L^0_s(\omega_2)}\right|\le
\sup_{s\le L_t^0(\omega_2)}\left|\frac{X_{\lfloor ns\rfloor}}{\sqrt n}-
W_{s}\right|,
\]
whence \eqref{eq:probx1}. Thus, putting
\[A_n:=\left\{(\omega_1,\omega_2):
\sup_{s\le t}\left|\frac{X_{\lfloor nL^0_s(\omega_2)\rfloor}
(\omega_1)}{\sqrt n}-
W_{L^0_s(\omega_2)}(\omega_1)\right|>\eps\right\},\]
we have that $\Prob_1(A_n)\to0$, by Fubini and the dominated convergence
theorem we get $\Prob_1\otimes\Prob_2(A_n)\to0$ and we are done.
\end{proof}

%%%%%%%%%%%%%%%%%%%%%%%%%%%%%%%%%%%%%%%%%
%%%	REMARKS     		      %%%
%%%%%%%%%%%%%%%%%%%%%%%%%%%%%%%%%%%%%%%%%

\section{Final remarks}
\label{sec:fr}

The results of the previous sections show that the random
walker on $\mathbf C_2$ spends most of the time walking
along the vertical direction.
One would ask to what extent the resemblance between the
simple random walk on $\Z$ and the vertical component
of the simple random walk on $\mathbf C_2$ is apparent.
The answer is that, to leading order, the expected values of the
distances after $n$ steps of these two walks are 
indistinguishable.
Indeed if we denote by $X_n$ the position at time $n$ of the
walker on $\Z$, one could easily compute
\begin{equation*}
\E[|X_n|]\stackrel{n\to\infty}\sim \sqrt{\frac\pi2}\,n^{1/2},
\end{equation*}
which coincides with the estimate of Theorem~\ref{th:disty}.
Moreover, a comparison between Theorem~\ref{th:mdy}
and \cite[Theorem 2.14]{cf:PP} and
between Theorem~\ref{th:msy} and \cite[Theorem 3.4]{cf:PP}
shows the same coincidence for the estimates of the
maximal deviation and the  maximal span respectively.

The inhomogeneity of $\mathbf C_2$ results in the difference
between the behaviour of the horizontal and vertical
components of the walk. Indeed while for any positive integer
$d$ the expected value of the distance from the origin of
the simple random walker on $\Z^d$,
after $n$ steps, is of order $n^{1/2}$ 
(this may be easily proven via a conditioning argument,
with respect to the proportions of time spent along each
of the $d$ main directions), the expected value of the horizontal
distance on $\mathbf C_2$ is of order $n^{1/4}$, that is on
this direction we have a subdiffusive behaviour. Of course
this is due to the delay observed on the $x$-axis
while the walker explores the teeth of the comb (recall that
the simple random walk on $\Z$, although recurrent, is 
zero-recurrent, that is the expected value of the first
return time to the origin is infinite).

The difference between the two projections of the walk is
remarkable also when one properly rescales the process as
we did in Section~\ref{sec:scaling}. Indeed 
the space where the ``rescaled walker'' lives is $\R^2$,
endowed with the topology for which any
path between two points $(x_1,y_1)$ and $(x_2,y_2)$ must
necessarily include the three edges $(x_1,y_1)$--$(x_1,0)$, 
$(x_1,0)$--$(x_2,0)$
and $(x_2,0)$--$(x_2,y_2)$. This reflects on the limiting process: 
the horizontal component may change only when the vertical one passes
through 0. In fact, as we proved in Theorem \ref{th:scaling},
$S_{nt}^y/n^{1/2}$ converges to a standard Brownian motion $B$, 
(the walk on the vertical direction is, to leading order, unaffected
by the bias of the random holding time at zero), while
$S_{nt}^x/n^{1/4}$ converges to a Brownian motion whose clock is the local
time at zero of $B$.

The fact that the walker on $\mathbf C_2$ essentially behaves
like the walker on $\Z$ makes it clear that $\delta_w(\mathbf C_2)$
must be 2, as we proved by combinatorial methods in 
Section~\ref{sec:maxdev} (note that $\delta_w(\Z^d)=2$ for all integer
$d\ge1$). Indeed it is very likely that the walker will exit the ball
of radius $n$ moving along some tooth of $\mathbf C_2$. This disproves
the Einstein relation between $\delta_s(\mathbf C_2)=3/2$, 
$\delta_f(\mathbf C_2)=2$ and $\delta_w(\mathbf C_2)=2$.
The failure of this relation in this case is due to the inhomogeneity
of this particular graph. Indeed for strongly inhomogeneous graphs
the growth exponent does not give an accurate description of the ``fractal
properties'' of the graph (the assignment of the same
$\delta_f$ of $\Z^2$ to $\mathbf C_2$ disregards the topology
of the two structures). It seems that defining $\delta_f$ as the
growth exponent of the graph makes sense only for homogeneous graphs
(like $\Z^d$) or self-similar graphs, which have a clear fractal nature
(like the Sierpi\'nski graph). It is our opinion that the study
of the sense in which a graph has a fractal nature and what is 
the proper definition of its fractal dimension  should 
require further investigations.

\section*{Acknowledgments}

I feel particularly indebted to Peter Grabner
and Helmut Prodinger for raising the questions discussed
in this paper and for the stimulating discussions during which
they suggested the techniques used in Sections~\ref{sec:distances}, 
\ref{sec:maxdev} and \ref{sec:maxspan}. 
I would also like to thank  Jean-Francois Le Gall for suggesting the limit of 
the rescaled process and giving me some precious hints.

\bigskip
%%%%%%%%%%%%%%%%%%%%%%%%%%%%%%%%%%%%%%%%%
%%%	BIBLIOGRAPHY		      %%%
%%%%%%%%%%%%%%%%%%%%%%%%%%%%%%%%%%%%%%%%%

\end{document}